\documentclass[12pt]{article}
\usepackage{amsfonts} 
\usepackage{latexsym} 

\def\ourstretch{1.3}
\renewcommand{\baselinestretch}{\ourstretch}

\newtheorem{thm}{Theorem}
\newtheorem{lem}{Lemma}

\newtheorem{df}{Definition}
\makeatletter
\newenvironment{pf}[1][]%
 {\def\proof@temp{#1}\par\noindent
  \textsc{Proof}\ifx\proof@temp\@empty\else\ (#1)\fi\hspace{1em}}
 {~~ \hphantom{mm}\hfill~~{$\Box$}\par\vspace{.2\baselineskip}}
\makeatother

\newcommand{\g}{\mathfrak{g}}
\newcommand{\h}{\mathfrak{h}}
\renewcommand{\k}{\mathfrak{k}}

\newcommand{\p}{\mathfrak{p}}
\renewcommand{\r}{\mathfrak{r}}
\newcommand{\s}{\mathfrak{s}}
\renewcommand{\O}{\mathcal{O}}

\newcommand{\BF}{B$\Rightarrow$F}
\newcommand{\fie}{\varphi}
\makeatletter
\newcommand{\Ad}{\mathop{\operator@font Ad}\nolimits}
\newcommand{\ad}{\mathop{\operator@font ad}\nolimits}
\makeatother

\title{Compact Coadjoint Orbits}
\author{%
~\\\large John Rawnsley\\[30pt]
Mathematics Institute\\
University of Warwick\\
Coventry\ \ CV4\ \ 7AL\\
United Kingdom\\[20pt]{\small Email: 
{\tt J.Rawnsley@warwick.ac.uk}}\\[20pt]
23 June, 2003
}
\date{
}

\begin{document}
\renewcommand{\baselinestretch}{1}
\maketitle
\begin{abstract}
I give an answer to the question ``Which groups have compact coadjoint
orbits?''. Whilst I thought that the answer, which is straightforward,
must be in the literature, I was unable to find it. This note aims to
rectify this.

It is also a plea: If the result is already published then I would like
to be told the reference.
\end{abstract}
\thispagestyle{empty}
%
\renewcommand{\baselinestretch}{\ourstretch}

\setcounter{page}{1}

\section{Introduction}

This note is a response to the question of ``Which connected Lie groups
have compact coadjoint orbits?'' which I was asked whilst on a visit to
the Differential Geometry group at the ULB in Brussels. The obvious
answer is that they are the coadjoint orbits of a compact Lie group. But
that is not completely true as the origin in $\g^*$ is a compact
coadjoint orbit in every Lie algebra. More generally, an element of the
dual of a Lie algebra which vanishes on the derived algebra will be a
fixed point of the coadjoint action and so give an orbit consisting of a
single point.

We aim to show by elementary means that these two cases essentially
account for all compact orbits in the following sense: if $\O$ is a compact
coadjoint orbit for the group $G$ then there is a closed normal subgroup
$H$ of $G$ which fixes each point of $\O$, and such that $G/H$ is a
compact semisimple Lie group in such a way that $\O$ is the sum of a
$G$-fixed element of $\h^*$ and a coadjoint orbit of $G/H$ pulled back
to $\g^*$.

I was surprised that I could not find a reference for such a basic
result despite the length of time for which coadjoint orbits have been
studied. This note is an attempt to give an elementary answer to the
question, and at the same time an appeal for references to the
literature in case  the answer is already known (although not to me, for
which I apologise in advance).

I would also like to thank the members of the group ``Mechanics,
Quantisation and Geometry'' for inviting me to Marseille, 
and helping to improve upon the first version of this note. However, any
remaining mistakes or inaccuracies are entirely mine.

We study the problem by restricting elements of a coadjoint orbit to
subalgebras. Such a restriction will not in general be a single orbit,
but will be a union of orbits. In the resulting set of orbits for the
subalgebra compactness of the original orbit may be lost since the new
orbits need not, a priori, be closed. In stead we weaken compactness to
boundedness which is preserved under restriction. So we study the
apparently more general notion of bounded coadjoint orbits. At the last
step we see that bounded orbits are necessarily compact, so we recover
the desired situation.

\section{Two technical lemmas}

A class of Lie algebras where there are no unbounded
coadjoint orbits are the abelian Lie algebras since there the coadjoint
action is trivial so each orbit consists of a single point fixed under
the whole group.

\begin{df} \rm
Say that a connected Lie group has property \BF\ if 
each bounded coadjoint orbit consists of a single fixed point.
\end{df}

Obviously, abelian Lie groups have property \BF. 

\begin{lem}\label{lem-restr}
Let $\h$ be an ideal in $\g$ and $r \colon \g^* \to \h^*$ be the
restriction map. Suppose that $\h$ has property \BF\ and $\O$ is a
bounded orbit in $\g^*$. Then $r(\O)$ is a single point, a fixed point
for the action of $G$ on $\h^*$ as automorphisms of $\h$.
\end{lem}

\begin{pf}
Let $H \subset G$ be the corresponding Lie groups and $\O\subset \g^*$ a
bounded coadjoint orbit of $G$. $r(\O)$ is a bounded set in $\h^*$ which
breaks up into bounded $H$ orbits. Since $\h$ has property \BF, these
must then be fixed points for $H$, so we have the basic result:
\[
f \in \O \quad \Rightarrow \quad h.(f\vert_{\h}) =
f\vert_{\h}, \quad \forall h \in H.
\]
Since $\h$ is an ideal in $\g$, $h.(f\vert_{\h}) = (h.f)\vert_{\h}$ and so
$h.f - f$ vanishes on $\h$.

In particular $h.f - f$ vanishes on $[\h,\g]$ using again that $\h$ is
an ideal. Then $h_1.(h_2.f-f) = h_2.f-f$ for all $h_1$, $h_2$ in $H$. If
we set $\fie(h) = h.f-f$ then $\fie(h_1h_2) = \fie(h_1) + \fie(h_2)$ so
that $\fie \colon H \to \g^*$ is a homomorphism of Lie groups to the
additive group of $\g^*$. The image lies in $\O-f$ and so is bounded.
But the only bounded subgroup of a vector space is the trivial subgroup
$\{0\}$. Hence $\fie(h) = 0$ for all $h$ and thus $h.f = f$. It follows
that $f$ vanishes on $[\h,\g]$ and hence that $g.f-f$ vanishes on $\h$
for all $g \in G$. Thus $r(\O)$ is a single point. Since $\h$ is an
ideal, $(g.f)|_{\h} = g.(f|_{\h})$ and the rest of the statement
follows.
\end{pf}

\begin{lem}\label{lem-BF}
Let
\[
0 \to \h \to \g \to \k \to 0
\]
be a short exact sequence of Lie algebras. If $\h$ and $\k$ have
property \BF\ then so has $\g$.
\end{lem}

\begin{pf}
We know from the previous Lemma that $r(f)$ is a single $G$-fixed point 
in $\h^*$.

Pick $f_0 \in \O$ and consider the $G$ orbit of $f-f_0$. Since $f$,
$g.f$, $g.f_0$ all equal to $f_0$ on $\h$, all points in the orbit of
$f-f_0$ vanish on $\h$ and so project to a $K$-orbit on $\k^*$. Since
$g.(f-f_0)$ lies in the bounded set $\O-\O$, we obtain a bounded orbit
in $\k^*$ and since $\k$ has property \BF, this is a single $K$-fixed
point. Hence $g.(f-f_0)=f-f_0$ for all $g \in G$ since $G$ induces the
$K$-action.

But $f=g'.f_0$ so $gg'.f_0 - g.f_0 = g'.f_0 - f_0$. Thus defining
$\fie(g) = g.f_0 - f_0$, we obtain a homomorphism $\fie \colon G \to
\g^*$ whose image is bounded, and hence is zero. Thus $g.f_0 =f_0$
and hence $\O$ consists of a single fixed point. Thus $G$ has
property \BF.
\end{pf}

\section{Solvable Lie Algebras}

\begin{thm}\label{thm-solv}
If $\g$ is solvable then $\g$ has property \BF.
\end{thm}

\begin{pf}
A solvable Lie algebra can be built from successive extensions by
abelian Lie algebras starting from an abelian Lie algebra, and, as we
have already observed, abelian Lie algebras trivially have property \BF.
Thus all solvable Lie algebras have property \BF.
\end{pf}

\section{Semisimple Lie Algebras}

Let $\g$ be a real semisimple Lie algebra and $\O \subset \g^*$ a coadjoint
orbit. Then using the Killing form, $\O$ corresponds with a conjugacy
class $\Ad G (\xi)$ in $\g$ which is bounded if and only if $\O$ is.
$\g$ is a direct sum of simple ideals $\g_i$, and $\xi$ a sum of
elements $\xi_i\in\g_i$. $\O$ will be bounded if and only if all the
conjugacy classes of the $\xi_i$ are bounded. So in looking for bounded
coadjoint orbits of semisimple Lie algebras we can assume $\g$ is
simple. We claim that a simple Lie algebra only has a non-trivial
bounded conjugacy class when the Lie algebra is compact or equivalently
that non-compact real simple Lie algebras have property \BF.

\begin{lem}\label{lem:commutes} 
Let $\g$ be a real semisimple Lie algebra and $\xi\in\g$ have a bounded
conjugacy class. If $\eta\in\g$ has $\ad\eta$ semisimple with only real
eigenvalues then $[\xi,\eta]=0$.
\end{lem}

\begin{pf}
Under the above assumptions, $\g$ will decompose as a direct sum of
eigenspaces $\g_\lambda$ for $\ad\eta$. If $\xi=\sum_\lambda
\xi_\lambda$ is the corresponding decomposition of $\xi$ then
\[
e^{t\ad\eta}\xi = \sum_\lambda e^{t\lambda} \xi_\lambda.
\]
The right hand side can only be bounded if $\xi_\lambda=0$ whenever
$\lambda\ne0$. Hence $[\eta,\xi] =0$.
\end{pf}

\begin{thm}\label{thm:semisimple}
If $\g$ is a non-compact simple Lie algebra then $\g$ has property
\BF. In this case, the only fixed point is the origin.
\end{thm}

\begin{pf}
A non-compact simple Lie algebra is generated by elements $\eta$ for
which $\ad \eta$ is semisimple with all eigenvalues real. To see this
just take any Cartan decomposition $\g = \k + \p$. For $\g$ non-compact
$\p\ne0$. But $[\k,\p]\subset \p$ and $[\p,\p] \subset \k$ imply that
$\p + [\p,\p]$ is a non-trivial ideal in $\g$ so $\g=\p + [\p,\p]$ since
$\g$ is simple. All elements $\eta$ of $\p$ are $\ad$-semisimple with
real eigenvalues. If $\xi$ corresponds under the Killing form with an
element of a bounded coadjoint orbit then we have $[\xi,\g]=0$ by Lemma
\ref{lem:commutes} and hence $\xi=0$.
\end{pf}

\section{The General Case}

Let $\g$ be an arbitrary Lie algebra and $\r$ its solvable radical. Then
$\s = \g/\r$ is semi-simple. We can separate $\s$ into compact and
non-compact ideals $\s = \s_c + \s_n$ and let $\g_n$ be the ideal of
$\g$ projecting onto $\s_n$. Then $\r \subset \g_n$ with $\g_n/\r \cong
\s_n$; since both $\r$ and $\s_n$ have property \BF{}  so does  $\g_n$
by Lemma \ref{lem-BF}. The quotient $\g/\g_n \cong \s_c$ is compact
semisimple. We shall show that a bounded coadjoint orbit $\O$ is
essentially a coadjoint orbit of the compact semisimple Lie algebra
$\s_c$.

Let $\O \subset \g$ be a bounded orbit. Then the restrictions $r(f)$ to
$\g_n$ of elements $f$ of $\O$ are all equal to a single $G$-fixed point
in $\g_n^*$ by Lemma \ref{lem-restr}. Pick $f_0 \in \O$ and consider
$f-f_0 \in \g^*$. Its orbit $\O'$ will also be bounded and all the
points in $\O'$ vanish on $\g_n$. It follows that $\O'$ projects
isomorphically to an orbit in the dual of $\g/\g_n$ and so is an orbit
of a compact semisimple Lie group.

If we choose a Levi factor so that $\s$ becomes a subalgebra of $\g$ and
$\g = \r + \s$ is a semidirect product then $\g_n = \r + \s_n$ and $\g =
\g_n + \s_c$. $f_0$ restricted to $\g_n$ is independent of $f_0$ and is
a $G$ fixed element $f_1$ of $\g_n^*$ which we can extend by 0 to $\g$.
Then $\O = f_1 + \O'$ where elements of $\O' \subset \s_c$ are extended
to $\g$ as zero on $\g_n$. In summary:

\begin{thm}
If $\O$ is a bounded coadjoint orbit of a Lie algebra, then $\O$ is
closed, hence compact and has the form $\O= f_1 + \O_c$ where $\O_c$ is
(the pull-back of) a coadjoint orbit of the compact semisimple quotient
and $f_1$ is a $G$-fixed element of the dual of the non-compact part
$\g_n$.
\end{thm}

\end{document}